\documentclass[]{article}
\usepackage[T2A]{fontenc}
\usepackage[cp1251]{inputenc}
\usepackage{amssymb}
\usepackage{amsmath}
\usepackage{amsthm}

\title{\bf Subtraction Menger algebras}
\author{W. A. Dudek and V. S. Trokhimenko}
\date{}

\begin{document}
\sloppy \maketitle

\newtheorem{theorem}{Theorem}
\newtheorem{proposition}{Proposition}
\newtheorem{definition}{Definition}
\newtheorem{collolary}{Corollary}
\newcommand{\pr}{\mbox{pr}_1\,}
\newcommand{\spr}{\mbox{\scriptsize pr}_1\,}

\begin{abstract}
Abstract characterizations of Menger algebras of partial $n$-place
functions  defined on a set $A$ and closed under the set-theoretic
difference functions treatment as subsets of the Cartesian product
$A^{n+1}$ are given.
\end{abstract}

\medskip {\bf 1.}
Let $A^n$ be the $n$-th Cartesian product of a set $A$. Any
partial mapping from $A^n$ into $A$ is called a {\it partial
$n$-place function}. The set of all such mappings is denoted by
$\mathcal{F}(A^n,A)$. On $\mathcal{F}(A^n,A)$ we define the {\it
Menger superposition} (composition) of $n$-place functions
$\mathrm{O}\colon\,(f,g_1,\ldots,g_n)\mapsto f[g_1\ldots g_n]$ as
follows:
\begin{equation}\label{e-1}
  (\bar{a},c)\in f[g_1\ldots
g_n]\longleftrightarrow(\exists\bar{b})\Big((\bar{a},b_1)\in
g_1\wedge\ldots\wedge(\bar{a},b_n)\in g_n\wedge(\bar{b},c)\in
f\Big)
\end{equation}
for all $\bar{a}\in A^n$, $\bar{b}=(b_1,\ldots,b_n)\in A^n$, $c\in
A$.

Each subalgebra $(\Phi,\mathrm{O})$, where
$\Phi\subset\mathcal{F}(A^n,A)$, of the algebra
$(\mathcal{F}(A^n,A),\mathrm{O})$ is a Menger algebra of rank $n$
in the sense of \cite{Dudtro2, Dudtro3, SchTro}. Menger algebras
of partial $n$-place functions are partially ordered by the
set-theoretic inclusion, i.e., such algebras can be considered as
algebras of the form $(\Phi,\mathrm{O},\subset)$. The first
abstract characterization of such algebras was given in
\cite{tro1}. Later, in \cite{tro2, TroSch} there have been found
abstract characterizations of Menger algebras of $n$-place
functions closed with respect to the set-theoretic intersection
and union of functions, i.e., Menger algebras of the form
$(\Phi,\mathrm{O},\cap)$, $(\Phi,\mathrm{O},\cup)$ and
$(\Phi,\mathrm{O},\cap,\cup)$.

As is well known, the set-theoretic inclusion $\subset$ and the
operations $\cap$, $\cup$ can be expressed by the set-theoretic
difference (subtraction) in the following way:
\begin{eqnarray}
& A\subset B\longleftrightarrow A\!\setminus B=\varnothing,\quad
A\cap B=A\!\setminus(A\!\setminus B),&\nonumber\\[4pt]
& A\cup B=C\!\setminus((C\!\setminus A)\cap(C\!\setminus
B)),&\nonumber
\end{eqnarray}
where $A,B,C$ are arbitrary sets such that $A\subset C$ and
$B\subset C$.

Thus it make sense to examine sets of functions closed with
respect to the subtraction of functions. Such sets of functions
are called {\it difference semigroups}, their abstract analogs --
{\it subtraction semigroups}. Properties of subtraction semigroups
were found in \cite{Abbott}. The investigation of difference
semigroups was initiated by B. M. Schein in \cite{Schein}.

Below we present a generalization of Schein's results to the case
of Menger algebras of $n$-place functions, i.e., to the case of
algebras $(\Phi,\mathrm{O},\!\setminus,\varnothing)$, where
$\Phi\subset\mathcal{F}(A^n, A)$, $\varnothing\in\Phi $. Such
algebras will be called \textit{difference Menger algebras}.

\medskip

{\bf 2.} A {\it Menger algebra of rank $n$} is a non-empty set $G$
with one $(n+1)$-ary operation $o(x,y_1,\ldots,y_n)=x[y_1\ldots
y_n]$ satisfying the identity:
\begin{equation}\label{e-2} x [y_1
\ldots y_n] [z_1\ldots z_n] = x [y_1 [z_1\ldots z_n]\ldots y_n
[z_1\ldots z_n]].
\end{equation}

A Menger algebra of rank $1$ is a semigroup. A Menger algebra
$(G,o)$ of rank $n$ is called {\it unitary}, if it contains {\it
selectors}, i.e., elements $e_1,\ldots,e_n\in G$, such that
$x[e_1\ldots e_n] = x$ and $e_i[x_1\ldots x_n]=x_i$ for all
$x,x_1,\ldots,x_n\in G$, $i=1,\ldots,n$. One can prove (see
\cite{Dudtro2, Dudtro3}), that every Menger algebra $(G,o)$ of
rank $n$ can be isomorphically embedded into the unitary Menger
algebra $(G^*,o^*)$ of the same rank with selectors
$e_1,\ldots,e_n\not\in G$ such that $G\cup\{e_1,\ldots,e_n\}$ is
the generating set of $(G^*,o^*)$.

Let $(G,o)$ be a Menger algebra of rank $n$. Let's consider the
alphabet $ G\cup\{[\,,\,],x\}$, where $[\,,\,],x$ does not belong
to $G$, and construct over this alphabet the set $T_n(G)$ of {\it
polynomials} such that:
\begin{itemize}
\item[$a)$] \ $x\in T_n(G)$;
\item[$b)$] \ if $i\in\{1,\ldots,n\}$, $a,b_1,\ldots,b_{i-1},b_{i+1},\ldots,b_n\in G$,\\
$t\in T_n(G)$, then $a[b_1\ldots b_{i-1}t\,b_{i+1},\ldots b_n]\in
T_n(G)$;
\item[$c)$] \ $T_n(G)$ contains those and only those polynomials which
are constructed by $a$) and $b)$.
\end{itemize}

A binary relation $\rho\subset G\times G$, where $(G,o)$ is a
Menger algebra of rank $n$, is
\begin{itemize}
\item {\it stable} if for all $x,y,x_i,y_i\in G$, $i=1,\ldots,n$
$$
(x,y),(x_1,y_1),\ldots,(x_n,y_n)\in\rho\longrightarrow
(x[x_1\ldots x_n], y[y_1\ldots y_n])\in\rho;
$$
\item {\it $l$-regular}, if for any $x,y,z_i\in G$, $i=1,\ldots,n$
$$
(x,y)\in\rho\longrightarrow (x[z_1\ldots z_n], y[z_1\ldots
z_n])\in\rho;
$$
\item {\it $v$-regular}, if for all $x_i,y_i,z\in G$, $i=1,\ldots,n$
$$
(x_1,y_1),\ldots,(x_n,y_n)\in\rho\longrightarrow (z[x_1\ldots x_n],z[y_1\ldots y_n])\in\rho;
$$
\item {\it $i$-regular} $(1\leqslant i\leqslant n)$, if for all $u,x,y\in G$, $\bar{w}\in G^n$
$$
(x,y)\in\rho\longrightarrow (u[\bar{w}|_ix],
u[\bar{w}|_iy])\in\rho;
$$
\item {\it weakly steady} if for all $x,y,z\in G$, $t_1,t_2\in T_n(G)$
$$
(x,y),(z,t_1(x)),(z,t_2(y))\in \rho\longrightarrow
(z,t_2(x))\in\rho,
$$
\end{itemize}
where $\bar{w}=(w_1,\ldots,w_n)$ and $u[\bar{w}|_i\,x]=u[w_1\ldots
w_{i-1}xw_{i+1}\ldots w_n]$. It is clear that a
quasiorder\footnote{\,Recall that a \textit{quasiorder} is a
reflexive and transitive binary relation.} on a Menger algebra is
$v$-regular if and only if it is $i$-regular for every
$i=1,\ldots,n$. A quasiorder is stable if and only if at the same
time it is $v$-regular and $l$-regular.

A subset $H$ of a Menger algebra $(G,o)$ is called
\begin{itemize}
\item {\it stable} if
$$
g,g_1,\ldots,g_n\in H\longrightarrow g[g_1\ldots g_n]\in H;
$$
\item an {\it $l$-ideal}, if for all $x,h_1,\ldots,h_n\in G$
$$
(h_1,\ldots,h_n)\in G^n\!\setminus (G\!\setminus H)^n
\longrightarrow x [h_1\ldots h_n]\in H;
$$
\item an {\it $i$-ideal} ($1\leqslant i\leqslant n$), if for all $h,u\in G$, $\bar{w}\in G^n $
$$
h\in H\longrightarrow u[\bar{w}|_i h]\in H.
$$
\end{itemize}
Clearly, $H$ is an $l$-ideal if and only if it is an $i$-ideal for
every $i=1,\ldots,n$.

\begin{definition}\label{D-sub}\rm An algebra $(G,-,0)$ of type $(2,0)$ is called a {\it subtraction algebra}
if it satisfies the following identities:
\begin{eqnarray}
 &&\label{e-3} x-(y-x)=x,\\
 &&\label{e-4}x-(x-y)=y-(y-x),\\
 &&\label{e-5}(x-y)-z=(x-z)-y,\\
 &&\label{e-6}0-0=0
\end{eqnarray}
for all $x,y,z\in G$.
\end{definition}

\begin{proposition} {\sc (Abbott \cite{Abbott})}\label{P-1}
Any subtraction algebra satisfies the identity:
\begin{equation}\label{e-7}
0=x-x.
\end{equation}
\end{proposition}
\begin{proof} Below we give a short proof of this identity:
\begin{eqnarray*}
0\!\!\!\!&\stackrel{(\ref{e-3})}{=}&\!\!\!\!
0-((0-(x-x))-0)\stackrel{(\ref{e-5})}{=}0-((0-0)-(x-x))\stackrel{(\ref{e-6})}{=}0-(0-(x-x))\\
 &\stackrel{(\ref{e-4})}{=}&\!\!\!\!(x-x)-((x-x)-0)\stackrel{(\ref{e-5})}{=}(x-x)-((x-0)-x)\\
 &\stackrel{(\ref{e-5})}{=}&\!\!\!\!(x-((x-0)-x))-x\stackrel{(\ref{e-3})}{=}x-x,
\end{eqnarray*}
was required to show.
\end{proof}

From \eqref{e-7}, by using \eqref{e-3}, we obtain the following
two identities:

\begin{equation}\label{e-8}
  x-0 = x,\quad 0-x = 0.
\end{equation}
Similarly, from \eqref{e-4}, \eqref{e-5}, \eqref{e-7} and
\eqref{e-8} we can deduce identities:
\begin{eqnarray}
&((x-y)-(x-z))-(z-y)=0,\label{e-9}\\
&(x-(x-y))-y = 0.\label{e-10}
\end{eqnarray}
Thus, any subtraction algebra $(G,-,0)$ is an implicative
BCK-algebra (cf. \cite{Iseki} or \cite{Kim}).

\begin{definition}\label{D-subalgMeng}\rm
An algebra $(G,o,-,0)$ of type $(n+1,2,0)$ is called a {\it
subtraction Menger algebra} of rank $n$, if $(G,o)$ is a Menger
algebra of rank $n,$ $(G,-,0) $ is a subtraction algebra and the
following conditions:
\begin{eqnarray}
 &&\label{e-11}(x-y)[z_1\ldots z_n]=x[z_1\ldots z_n]-y[z_1\ldots z_n],\\
 &&\label{e-12}u[\bar{w}|_i\,(x-(x-y))]=u[\bar{w}|_i\,x]-u[\bar{w}|_i\,x-y],\\
 &&\label{e-13}x-y=0\wedge z-t_1(x)=0\wedge z-t_2(y)=0\longrightarrow z-t_2(x)=0
\end{eqnarray}
are satisfied  for all $x,y,z,u,z_1,\ldots,z_n\in G$, $\bar{w}\in
G^n$, $i=1,\ldots,n$ and $t_1,t_2\in T_n(G)$.
\end{definition}

By putting $n=1$ in the above definition we obtain a \textit{weak
subtraction semigroup}\,\footnote{A weak subtraction semigroup
$(S,\cdot,-)$ is a semigroup $(S,\cdot)$ satisfying the identities
$(\ref{e-3})$, $(\ref{e-4})$, $(\ref{e-5})$, $x(y-z)=xy-xz$ and
$(x-(x-y))z=xz-(x-y)z$.} studied by B. M. Schein (cf.
\cite{Schein}). Such semigroups are isomorphic to some subtraction
semigroups of the form $(\Phi,\circ,\!\setminus)$.

\medskip

{\bf 3.} Now we can present the first result of our
paper.

\begin{theorem}\label{T-1}
Each difference Menger algebra of $n$-place functions is a subtraction Menger algebra
of rank $n$.
\end{theorem}
\begin{proof}
Let $(\Phi,\mathrm{O},\!\setminus,\varnothing)$ be a difference
Menger algebra of  $n$-place functions defined on $A$. Since, as
it is proved in \cite{Dudtro2}, the superposition $\mathrm{O}$
satisfies \eqref{e-2}, the algebra $(\Phi,\mathrm{O})$ is a Menger
algebra of rank $n$. From the results proved in \cite{Abbott} it
follows that the operation $\!\setminus$ satisfies \eqref{e-3},
\eqref{e-4} and \eqref{e-5}. Hence
$(\Phi,\!\setminus,\varnothing)$ is a subtraction algebra. Thus,
$(\Phi,\mathrm{O},\!\setminus,\varnothing)$ will be a subtraction
Menger algebra if \eqref{e-11}, \eqref{e-12} and \eqref{e-13} will
be satisfied.

To verify \eqref{e-11} observe that for each $(\bar{a},c)\in
(f\!\setminus g)[h_1\ldots h_n]$, where
$f,g,h_1,\ldots,h_n\in\Phi$, $\bar{a}\in A^n$, $c\in A$ there
exists $\bar{b}=(b_1,\ldots,b_n)\in A^n$ such that $(\bar{b},c)\in
f\!\setminus g$ and $(\bar{a},b_i)\in h_i$ for each
$i=1,\ldots,n$. Consequently, $(\bar{b},c)\in f$ and
$(\bar{b},c)\not\in g$. Thus, $(\bar{a},c)\in f[h_1\ldots h_n]$.
If $(\bar{a},c)\in g[h_1\ldots h_n]$, then there exists
$\bar{d}=(d_1,\ldots, d_n)\in A^n$ such that $(\bar{d},c)\in g$
and $(\bar{a},d_i)\in h_i$ for every $i=1,\ldots,n$. Since
$h_1,\ldots,h_n$ are functions, we obtain $b_i=d_i$ for all
$i=1,\ldots,n$. Thus $\bar{b}=\bar{d}$. Therefore $(\bar{b},c)\in
g$, which is impossible. Hence $(\bar{a},c)\not\in g[h_1\ldots
h_n]$. This means that $(\bar{a},c)\in f[h_1\ldots h_n]\setminus
g[h_1\ldots h_n]$. So, the following implication
\[
(\bar{a},c)\in(f\!\setminus g)[h_1\ldots h_n]\longrightarrow
(\bar{a},c)\in f[h_1\ldots h_n]\setminus g[h_1\ldots h_n]
\]
is valid for any $\bar{a}\in A^n$, $c\in A$, i.e., $(f\!\setminus
g)[h_1\ldots h_n]\subset f[h_1\ldots h_n]\setminus g[h_1\ldots
h_n]$.

Conversely, let $(\bar{a},c)\in f[h_1\ldots h_n]\setminus
g[h_1\ldots h_n]$. Then $(\bar{a},c)\in f[h_1\ldots h_n]$ and
$(\bar{a},c)\not\in g [h_1\ldots h_n]$. Thus, there exists
$\bar{b}=(b_1,\ldots,b_n)\in A^n$ such that $(\bar{b},c)\in f$,
$(\bar{b},c)\not\in g$ and $(\bar{a},b_i)\in h_i$ for each
$i=1,\ldots,n$. Hence, $(\bar{b},c)\in f\!\setminus g$ and
$(\bar{a},c)\in (f\!\setminus g)[h_1\ldots h_n]$. So,
\[
(\bar{a},c)\in f[h_1\ldots h_n]\setminus g[h_1\ldots
h_n]\longrightarrow (\bar{a},c)\in (f\!\setminus g)[h_1\ldots h_n]
\]
for any $\bar{a}\in A^n$, $c\in A$, i.e., $f[h_1\ldots
h_n]\setminus g[h_1\ldots h_n]\subset (f\!\setminus g)[h_1\ldots
h_n]$. Thus,
$$
(f\!\setminus g)[h_1\ldots h_n]=f[h_1\ldots h_n]\setminus
g[h_1\ldots h_n],
$$
which proves \eqref{e-11}.

Now, let $(\bar{a},c)\in u[\bar{\omega}|_i(f\!\setminus
(f\!\setminus g))] =u[\bar{\omega}|_i(f\cap g)]$, where
$f,g,u\in\Phi$, $\bar{\omega}\in\Phi^n$, $\bar{a}\in A^n$, $c\in
A$. Then there exists $\bar{b}=(b_1,\ldots,b_n)\in A^n$ such that
$(\bar{a},b_i)\in f\cap g$, $(\bar{a},b_j)\in\omega_j$,
$j\in\{1,\ldots,n\}\!\setminus\{i\}$ and $(\bar{b},c)\in u$. Since
$(\bar{a},b_i)\in f\cap g$ implies $(\bar{a},b_i)\not\in
f\!\setminus g$, we have $(\bar{a},c)\in u[\bar{\omega}|_if]$ and
$(\bar{a},c)\not\in u[\bar{\omega}|_i(f\!\setminus g)]$. Therefore
$(\bar{a},c)\in u[\bar{\omega}|_if]\setminus
u[\bar{\omega}|_i(f\!\setminus g)]$. Thus, we have shown that for
any $\bar{a}\in A^n$, $c\in A$ holds the implication
\[
(\bar{a},c)\in u[\bar{\omega}|_i(f\!\setminus (f\!\setminus
g))]\longrightarrow (\bar{a},c)\in u[\bar{\omega}|_if]\setminus
u[\bar{\omega}|_i(f\!\setminus g)],
\]
which is equivalent to the inclusion
$u[\bar{\omega}|_i(f\!\setminus (f\!\setminus g))]\subset
u[\bar{\omega}|_if]\setminus u[\bar{\omega}|_i(f\!\setminus g)]$.

Conversely, let $(\bar{a},c)\in u[\bar{\omega}|_if]\setminus
u[\bar{\omega}|_i(f\!\setminus g)]$. Then $(\bar{a},c)\in
u[\bar{\omega}|_if]$ and $(\bar{a},c)\not\in
u[\bar{\omega}|_i(f\!\setminus g)]$. The first of these two
conditions means that there exists $\bar{b}=(b_1,\ldots,b_n)\in
A^n$ such that $(\bar{a},b_i)\in f$, $(\bar{a},b_j)\in\omega_j$
for each $j\in\{1,\ldots,n\}\!\setminus\{i\}$ and $(\bar{b},c)\in
u$. It is easy to see that the second condition
$(\bar{a},c)\not\in u[\bar{\omega}|_i(f\!\setminus g)]$ is
equivalent to the implication
\begin{equation}\label{e-14}
  (\forall\bar{d})\Big((\bar{a},d_i)\in f\wedge\bigwedge_{j=1,j\neq i}^{n}(\bar{a},d_j)\in\omega_j
 \wedge(\bar{d},c)\in u\longrightarrow(\bar{a},d_i)\in g\Big),
\end{equation}
where $\bar{d}=(d_1,\ldots,d_n)\in A^n$. From this implication for
$\bar{d}=\bar{b}$, we obtain
\[
(\bar{a},b_i)\in f\wedge\bigwedge_{j=1,j\neq
i}^n(\bar{a},b_j)\in\omega_j\wedge (\bar{b},c)\in u\longrightarrow
(\bar{a},b_i)\in g,
\]
which gives $(\bar{a},b_i)\in g$. Therefore $(\bar{a},b_i)\in
f\cap g=f\!\setminus (f\!\setminus g)$. This means that
$(\bar{a},c)\in u[\bar{\omega}|_i(f\!\setminus (f\!\setminus
g))]$. So, the implication
\[
(\bar{a},c)\in u[\bar{\omega}|_if]\setminus
u[\bar{\omega}|_i(f\!\setminus g)]\longrightarrow (\bar{a},c)\in
u[\bar{\omega}|_i(f\!\setminus (f\!\setminus g))]
\]
is valid for all $\bar{a}\in A^n$, $c\in A$. Hence
$u[\bar{\omega}|_if]\setminus u[\bar{\omega}|_i(f \!\setminus
g)]\subset u[\bar{\omega}|_i(f\!\setminus (f\!\setminus g))]$.
Thus
\[
u[\bar{\omega}|_i(f\!\setminus (f \!\setminus g))]=
u[\bar{\omega}|_if]\setminus u[\bar{\omega}|_i(f\!\setminus g)].
\]
This proves \eqref{e-12}.

To prove \eqref{e-13} suppose that for some $f,g,h\in\Phi$ and
$t_1,t_2\in T_n(\Phi)$ we have $f\!\setminus g =\varnothing$,
$h\!\setminus t_1(f)=\varnothing$ and $h\!\setminus t_2(g)=
\varnothing$. Then $f\subset g$, $h\subset t_1(f)$ and $h\subset
t_2(g)$. Hence $f=g\circ\bigtriangleup_{\spr f}$ and $\pr
h\subset\pr f$, where $\pr f$ denotes the domain of $f$ and
$\Delta_{\spr f}$ is the identity binary relation on $\pr f$.

From the inclusion $h\subset t_2(g)$ we obtain
$$
h=h\circ\Delta_{\spr f}\subset t_2(g)\circ\Delta_{\spr
f}=t_2(g\circ\Delta_{\spr f})=t_2(f),
$$
which means that \eqref{e-13} is also satisfied. This completes
the proof that $(\Phi,\mathrm{O},\setminus,\varnothing)$ is a
subtraction Menger algebra of rank $n$.
\end{proof}

To prove the converse statement, we need to consider a number of
properties of a subtraction Menger algebra of rank $n$, introduce
some definitions and prove some auxiliary propositions.

\medskip {\bf 4.} Let $(G,o,-, 0)$ be a subtraction Menger
algebra of rank $n$.

\begin{proposition}\label{P-2}
In any subtraction Menger algebra of rank $n$ we have
$$
0[x_1\ldots x_n]=0, \ \ \ \ x[x_1\ldots x_{i-1}0\,x_{i+1}\ldots x_n]=0
$$
for all $x,x_1,\ldots,x_n\in G$, $i=1,\ldots,n$.
\end{proposition}
\begin{proof} Indeed, using \eqref{e-7} and \eqref{e-11} we obtain
$$
0[x_1\ldots x_n]=(0-0)[x_1\ldots x_n]=0[x_1\ldots
x_n]-0[x_1\ldots x_n]=0.
$$
Similarly, applying \eqref{e-12} and \eqref{e-7} we get
$$
u[\bar{w}|_i\,0]=u[\bar{w}|_i(0-(0-0))]=u[\bar{w}|_i\,0]-u[\bar{w}|_i(0-0)]=u[\bar{w}|_i\,0]-
u[\bar{w}|_i\,0]=0,
$$
which was to show.
\end{proof}

Let $\omega$ be a binary relation defined on $(G,o,-, 0)$ in the
following way:
\[
\omega =\{(x, y)\in G\times G\,|\,x-y = 0\}.
\]
Using \eqref{e-7}, \eqref{e-8} and \eqref{e-9} it is easy to see
that this is an order, i.e., a reflexive, transitive and
antisymmetric relation. In connection with this fact we will
sometimes write $x\leqslant y$ instead of $(x,y)\in\omega$. Using
this notation it is not difficult to verify that
\begin{eqnarray}
&&\label{e-15} 0\leqslant x,\quad x-y\leqslant x, \\
&&\label{e-16} x\leqslant y\longleftrightarrow x-(x-y) = x,\\
&&\label{e-17} x\leqslant y\longrightarrow x-z\leqslant y-z,\\
 &&\label{e-18} x\leqslant y\longrightarrow z-y\leqslant z-x,\\
 &&\label{e-19} x\leqslant y\wedge u\leqslant v\longrightarrow x-v\leqslant y-u
\end{eqnarray}
holds for all $x,y,z,u,v\in G$.

Moreover, in a subtraction algebra the following two identities
\begin{eqnarray}
 &&\label{e-20} (x-y)-y = x-y,\\
 &&\label{e-21} (x-y)-z = (x-z) - (y-z)
\end{eqnarray}
are valid (cf. \cite{Abbott, Iseki, Kim}).

\begin{proposition}\label{P-3} On the algebra $(G,o,-,0)$ the relation $\omega $
is stable and weakly steady.
\end{proposition}
\begin{proof}
Let $x\leqslant y$ for some $x,y\in G$. Then $x-y=0$ and
\[
(x-y)[z_1\ldots z_n]=0[z_1\ldots z_n]=(0-0)[z_1\ldots z_n]=0[z_1\ldots z_n]-0[z_1\ldots z_n]=0
\]
for all $z_1,\ldots,z_n\in G$. This, by \eqref{e-11}, implies
\[
x[z_1\ldots z_n]-y[z_1\ldots z_n]=0,
\]
i.e., $x[z_1\ldots z_n]\leqslant y[z_1\ldots z_n]$. Thus, $\omega$
is $l$-regular.

Moreover, from $x\leqslant y$, using \eqref{e-8}, we obtain
$x-(x-y)=x$, which together with \eqref{e-4}, gives $y-(y-x)=x$.
Consequently, for any $ u\in G $, $\bar{w}\in G^n $ we have
$u[\bar{w}|_i(y-(y-x))]=u[\bar{w}|_i\,x]$. This and \eqref{e-11}
give $u[\bar{w}|_i\,y]-u[\bar{w}|_i(y-x)]=u[\bar{w}|_i\,x]$.
Hence, according to \eqref{e-15}, we obtain
$u[\bar{w}|_i\,x]\leqslant u[\bar{w}|_i\,y]$. Thus, $\omega$ is
$i$-regular for every $i=1,\ldots,n$. Since $\omega$ is a
quasiorder, the last means that $\omega$ is $v$-regular. But
$\omega$ also is $l$-regular, hence it is stable.

It is clear that $\omega$ is weakly steady if and only if it
satisfies \eqref{e-13}.\footnote{\,In the case of semigroups the
fact that $\omega$ is weakly steady can be deduced directly from
the axioms of a weak subtraction semigroup (cf. \cite{Schein}).}
\end{proof}

\begin{proposition}\label{P-4} The
axiom $\eqref{e-12}$ is equivalent to each of the following
conditions:
\begin{eqnarray}
&&\label{e-22} x\leqslant y\longrightarrow u[\bar{w}|_i\,(y-x)]=
u[\bar{w}|_i\,y]-u[\bar{w}|_i\,x],\\
&&\label{e-23} x\leqslant y\longrightarrow t(y-x)=t(y)-t(x),\\
&&\label{e-24} t(x-(x-y))=t(x)-t(x-y)
\end{eqnarray} for all $x,y,u\in G$, $\bar{w}\in G^n$, $i=1,\ldots,n$, $t\in T_n(G)$.
\end{proposition}
\begin{proof} $\eqref{e-12}\longrightarrow\eqref{e-22}$. \ Suppose that the condition $\eqref{e-12}$ is satisfied and
$x\leqslant y$ for some $x,y\in G$. Then, according to
\eqref{e-16}, we have $x-(x-y)=x$. Hence, by \eqref{e-4}, we
obtain $y-(y-x)=x$. Thus, $y-x=y-(y-(y-x))$, which, in view of
\eqref{e-12}, gives
$u[\bar{w}|_i\,(y-x)]=u[\bar{w}|_i\,(y-(y-(y-x)))]=u[\bar{w}|_iy]-u[\bar{w}|_i(y-(y-x))]
=u[\bar{w}|_iy]-u[\bar{w}|_ix]$. This means that \eqref{e-12}
implies \eqref{e-22}.

$\eqref{e-22}\longrightarrow\eqref{e-23}$. \ From \eqref{e-22} it
follows that for $x\leqslant y$ and all polynomials $t\in T_n(G)$
of the form $t(x)=u[\bar{w}|_ix]$ the condition \eqref{e-23} is
satisfied. To prove that \eqref{e-23} is satisfied by an arbitrary
polynomial from $T_n(G)$ suppose that it is satisfied by some
$t'\in T_n(G)$. Since the relation $\omega$ is stable on the
algebra $(G,o,-,0)$, from $x\leqslant y$ it follows
$t'(x)\leqslant t'(y)$, which in view of \eqref{e-22}, implies
$$
u[\bar{w}|_i\,(t'(y)-t'(x))]=u[\bar{w}|_i\,t'(y)]-u[\bar{w}|_i\,t'(x)].
$$
But according to the assumption on $t'$ for $x\leqslant y$ we have
$t'(y)-t'(x)=t'(y-x)$, so the above equation can be written as
$$
u[\bar{w}|_i\,t'(y-x)]=u[\bar{w}|_i\,t'(y)]-u[\bar{w}|_i\,t'(x)].
$$
Thus, \eqref{e-23} is satisfied by polynomials of the form
$t(x)=u[\bar{w}|_it'(x)]$.

From the construction of $T_n(G)$ it follows that \eqref{e-23} is
satisfied by all polynomials $t\in T_n(G)$. Therefore \eqref{e-22}
implies \eqref{e-23}.

$\eqref{e-23}\longrightarrow\eqref{e-24}$. \ Since, by
\eqref{e-15}, $x-y\leqslant x$ holds for all $x,y\in G$, from
\eqref{e-23} it follows $t(x-(x-y))=t(x)-t(x-y)$ for any
polynomial $t\in T_n(G)$. Thus, \eqref{e-23} implies \eqref{e-24}.

$\eqref{e-24}\longrightarrow\eqref{e-12}$. \ By putting
$t(x)=u[\bar{w}|_i\,x]$ we obtain \eqref{e-12}.
\end{proof}

On a subtraction Menger algebra $(G,o,-,0)$ of rank $n$ we can
define a binary operation $\curlywedge$ by putting:
\begin{equation}\label{e-25}
  x\curlywedge y\stackrel{def}{=}x-(x-y).
\end{equation}
By using this operation the conditions \eqref{e-11}, \eqref{e-16},
\eqref{e-24} can be written in a more useful form:
\begin{eqnarray}
&&\label{e-26} u[\bar{w}|_i(x\curlywedge y)]=u[\bar{w}|_i\,x]-u[\bar{w}|_i(x-y)],\\
&&\label{e-27} x\leqslant y\longleftrightarrow x\curlywedge y=x,\\
&&\label{e-28} t(x\curlywedge y)=t(x)-t(x-y),
\end{eqnarray}
where $x,y,u\in G$, $\bar{w}\in G^n$, $i=1,\ldots,n$, $t\in
T_n(G)$. Moreover, from \eqref{e-11} and \eqref{e-25}, we can
deduce the identity:
\begin{equation}\label{e-29}
 (x\curlywedge y)[z_1\ldots z_n]=
 x[z_1\ldots z_n]\curlywedge y[z_1\ldots z_n].
\end{equation}

The algebra $(G,\curlywedge)$ is a lower semilattice. Directly
from the conditions  \eqref{e-3}~--~\eqref{e-10} we obtain (cf.
\cite{Abbott}) the following properties:
\begin{eqnarray}
&&\label{e-30}
 x\leqslant y\wedge x\leqslant z\longrightarrow x\leqslant y\curlywedge z,\\
&&\label{e-31}
 x\leqslant y\longrightarrow x\curlywedge z\leqslant y\curlywedge z,\\
&&\label{e-32} x\curlywedge y=0\longrightarrow x-y=x,\\
&&\label{e-33} (x-y)\curlywedge y=0,\\
&&\label{e-34} x\curlywedge(y-z)=(x\curlywedge y)-(x\curlywedge z),\\
&&\label{e-35} x-y=x-(x\curlywedge y),\\
&&\label{e-36} (x\curlywedge y)-(y-z)=x\curlywedge y\curlywedge z,\\
&&\label{e-37} (x\curlywedge y)-z=(x-z)\curlywedge(y-z),\\
&&\label{e-38} (x\curlywedge y)-z=(x-z)\curlywedge y
\end{eqnarray}
for all $x,y,z\in G$.
\begin{proposition}\label{P-5}
In a subtraction Menger algebra $(G,o,-,0)$ of rank $n$ the
following conditions
\begin{eqnarray}
&&\label{e-39} t(x-y)=t(x)-t(x\curlywedge y),\\
&&\label{e-40} t(x)-t(y)\leqslant t(x-y)
\end{eqnarray}
are valid for each $t\in T_n(G)$ and $x,y\in G$.
\end{proposition}

\begin{proof}
From \eqref{e-35} we obtain $t(x-y)=t(x-(x\curlywedge y))$ for
every $t\in T_n(G)$. \eqref{e-25} and \eqref{e-15} imply
$x\curlywedge y\leqslant x$, which together with \eqref{e-23}
gives $t(x-(x\curlywedge y))=t(x)-t(x\curlywedge y)$. Hence,
$t(x-y)=t(x)-t(x\curlywedge y)$. This proves \eqref{e-39}.

Since $x\curlywedge y\leqslant y$, the stability of $\omega$
implies $t(x\curlywedge y)\leqslant t(y)$ for every $t\in T_n(G)$.
From this, by applying \eqref{e-15} and \eqref{e-18}, we obtain
$t(x)-t(y)\leqslant t(x)-t(x\curlywedge y)=t(x-y)$, which proves
\eqref{e-40}.
\end{proof}

By $[0,a]$ we denote the {\it initial segment} of the algebra
$(G,-,0)$, i.e., the set of all $x\in G$ such that $0\leqslant
x\leqslant a$. According to \cite{Schein}, on any $[0,a]$ we can
define a binary operation $\curlyvee$ by putting:
\begin{equation}\label{e-41}
  x\curlyvee y\stackrel{def}{=}a-((a-x)\curlywedge(a-y))
\end{equation}
for all $x,y\in [0,a]$. It is not difficult to see that this
operation is idempotent and commutative, and $0$ is its neutral
element, i.e., $x\curlyvee x=x$, $x\curlyvee y=y\curlyvee x$,
$x\curlyvee 0=x$ for all $x,y\in [0,a]$.
\begin{proposition}\label{P-6}
For any $x,y\in[0,b]\subset[0,a]$, where $a,b\in G$, we have
\begin{equation}\label{e-42}
b-((b-x)\curlywedge(b-y))=a-((a-x)\curlywedge(a-y)).
\end{equation}
\end{proposition}
\begin{proof}
Note first that $b=b\curlywedge a$ because $b\leqslant a$.
Moreover, from $x\leqslant b$ and $y\leqslant b$, according to
\eqref{e-18}, we obtain $a-b\leqslant a-x$ and $a-b\leqslant a-y$.
This together with \eqref{e-30} gives
$a-b\leqslant(a-x)\curlywedge(a-y)$. Thus,
$(a-b)-((a-x)\curlywedge(a-y))=0$.

By \eqref{e-15} we have $b-((a-x)\curlywedge(a-y))\leqslant b$,
which implies
\begin{equation}\label{e-43}
  b\curlywedge(b-((a-x)\curlywedge(a-y)))=b-((a-x)\curlywedge(a-y)).
\end{equation}
Obviously $b=b\curlywedge b=b\curlywedge a$, $x=b\curlywedge x$,
$y=b\curlywedge y$. Therefore:\footnote{\,To reduce the number of
brackets we will write $x\curlywedge y-z$ instead of
$(x\curlywedge y)-z$.}
\begin{eqnarray*}
&&\!\!\!\!  b-((b-x)\curlywedge(b-y))=b\curlywedge
b-((b\curlywedge a-b\curlywedge x)\curlywedge(b\curlywedge
a-b\curlywedge y))\\
&\stackrel{\eqref{e-34}}{=}&\!\!\!\!
b\curlywedge b-(b\curlywedge(a-x)\curlywedge b\curlywedge(a-y))=
b\curlywedge b-b\curlywedge((a-x)\curlywedge(a-y))\\
&\stackrel{\eqref{e-34}}{=}&\!\!\!\!
b\curlywedge(b-((a-x)\curlywedge(a-y)))\stackrel{\eqref{e-42}}{=}
b-((a-x)\curlywedge(a-y))\\
&=&\!\!\!\! a\curlywedge b-((a-x)\curlywedge(a-y))
\stackrel{\eqref{e-25}}{=}(a-(a-b))-((a-x)\curlywedge(a-y))\\
&\stackrel{\eqref{e-21}}{=}&\!\!\!\!
(a-((a-x)\curlywedge(a-y)))-((a-b)-((a-x)\curlywedge(a-y)))\\
&=&\!\!\!\!(a-((a-x)\curlywedge(a-y)))-0\stackrel{\eqref{e-8}}{=}a-((a-x)\curlywedge(a-y)),
\end{eqnarray*}
which completes the proof.
\end{proof}
\begin{collolary}\label{C-1}
The condition $\eqref{e-42}$ is valid for all $x,y\in [0,a]\cap
[0,b]$.
\end{collolary}
\begin{proof}
Since $[0,a]\cap[0,b]=[0,a\curlywedge b]\subset[0,a]\cup[0,b]$, by
Proposition \ref{P-6}, for all $x,y\in [0,a]\cap[0,b]$ we have:
\begin{eqnarray}
&& a-((a-x)\curlywedge(a-y))=a\curlywedge b-((a\curlywedge
b-x)\curlywedge(a\curlywedge b-y)),\nonumber\\
&& b-((b-x)\curlywedge(b-y))=a\curlywedge b-((a\curlywedge
b-x)\curlywedge(a\curlywedge b-y)).\nonumber
\end{eqnarray}
This implies \eqref{e-42}.
\end{proof}

From the above corollary it follows that the value of $x\curlyvee
y$, if it exists, does not depend on the choice of the interval
$[0,a]$ containing the elements $x$ and $y$. In \cite{Abbott} it
is proved that for $x,y,z\in [0,a]$ we have:
\begin{eqnarray}
&&\label{e-44} x\curlywedge(x\curlyvee y)=x,\\
&&\label{e-45} x\curlyvee(x\curlywedge y)=x,\\
&&\label{e-46} (x\curlyvee y)\curlyvee z=x\curlyvee(y\curlyvee z),\\
&&\label{e-47} x\curlywedge(y\curlyvee z)=(x\curlywedge
y)\curlyvee(x\curlywedge z),\\
&&\label{e-48}
x\curlyvee(y\curlywedge z)=(x\curlyvee y)\curlywedge(x\curlyvee z),\\
&&\label{e-49} (x\curlyvee y)-z=(x-z)\curlyvee(y-z),\\
&&\label{e-50} x\leqslant z\wedge
y\leqslant z\longrightarrow x\curlyvee y\leqslant z,\\
&&\label{e-51} y\leqslant x\longrightarrow x=(x-y)\curlyvee y,\\
&&\label{e-52} x=(x\curlyvee y)-(y-x),\\
&&\label{e-53} x=(x\curlywedge y)\curlyvee(x-y).
\end{eqnarray}
From \eqref{e-44} it follows $x\leqslant x\curlyvee y$.

\begin{proposition}\label{P-7}
If for some $x,y\in G$ there exists $x\curlyvee y$, then for all
$u\in G$, $\bar{z},\bar{w}\in G^n$, $i=1,\ldots,n$ there are also
elements $x[\bar{z}]\curlyvee y[\bar{z}]$ and
$u[\bar{w}|_i\,x]\curlyvee u[\bar{w}|_i\,y]$, and the following
identities are satisfied:
\begin{eqnarray}
&&\label{e-55} (x\curlyvee y)[\bar{z}]=x[\bar{z}]\curlyvee
y[\bar{z}],\\
&&\label{e-56} u[\bar{w}|_i(x\curlyvee y)]=
u[\bar{w}|_i\,x]\curlyvee u[\bar{w}|_i\,y].
\end{eqnarray}
\end{proposition}
\begin{proof}
Suppose that an element $x\curlyvee y$ exists. Then $x\leqslant a$
and $y\leqslant a$ for some $a\in G$, which, by the $l$-regularity
of the relation $\omega$, implies $x[\bar{z}]\leqslant a[\bar{z}]$
and $y[\bar{z}]\leqslant a[\bar{z}]$ for any $\bar{z}\in G^n$.
This means that $x[\bar{z}]\curlyvee y[\bar{z}]$ exists and
\begin{eqnarray}
&& (x\curlyvee
y)[\bar{z}]\stackrel{\eqref{e-41}}{=}(a-((a-x)\curlywedge(a-y)))[\bar{z}]\stackrel{\eqref{e-11}}{=}
a[\bar{z}]-((a-x)\curlywedge(a-y))[\bar{z}]\nonumber\\
 &&\stackrel{\eqref{e-29}}{=}a[\bar{z}]-((a-x)[\bar{z}]\curlywedge(a-y)[\bar{z}])\stackrel{\eqref{e-11}}{=}
a[\bar{z}]-((a[\bar{z}]-x[\bar{z}])\curlywedge(a[\bar{z}]-y[\bar{z}]))\nonumber\\
 &&\stackrel{\eqref{e-41}}{=}x[\bar{z}]\curlyvee y[\bar{z}].\nonumber
\end{eqnarray}
This proves \eqref{e-55}.

Further, from $x\leqslant a$, $y\leqslant a$ and the
$i$-regularity of $\omega$ we obtain $u[\bar{w}|_i\,x]\leqslant
u[\bar{w}|_i\,a]$ and $u[\bar{w}|_i\,y]\leqslant
u[\bar{w}|_i\,a]$. Hence, there exists an element
$u[\bar{w}|_i\,x]\curlyvee u[\bar{w}|_i\,y]$. Since $x\leqslant
x\curlyvee y$ and $y\leqslant x\curlyvee y$, we also have
$u[\bar{w}|_i\,x]\leqslant u[\bar{w}|_i(x\curlyvee y)]$ and
$u[\bar{w}|_i\,y]\leqslant u[\bar{w}|_i(x\curlyvee y)]$, which,
according to \eqref{e-50}, gives
\begin{equation}\label{e-57}
u[\bar{w}|_i\,x]\curlyvee u[\bar{w}|_i\,y]\leqslant
u[\bar{w}|_i(x\curlyvee y)].
\end{equation}

On the other side, the existence of $u[\bar{w}|_i\,x]\curlyvee
u[\bar{w}|_i\,y]$ implies,
$$
u[\bar{w}|_i\,x]\leqslant u[\bar{w}|_i\,x]\curlyvee
u[\bar{w}|_i\,y] \ \mbox{ and }\ u[\bar{w}|_i\,y]\leqslant
u[\bar{w}|_i\,x]\curlyvee u[\bar{w}|_i\,y].
$$
Moreover,
\[
u[\bar{w}|_i(x\curlyvee
y)]-u[\bar{w}|_i(y-x)]\stackrel{\eqref{e-40}}{\leqslant}
u[\bar{w}|_i((x\curlyvee y)-(y-x))]
\stackrel{\eqref{e-52}}{=}u[\bar{w}|_i\,x].
\]
Consequently,
\begin{equation}\label{e-58}
u[\bar{w}|_i(x\curlyvee y)]-u[\bar{w}|_i(y-x)]\leqslant
u[\bar{w}|_i\,x]\curlyvee u[\bar{w}|_i\,y].
\end{equation}

But $y-x\leqslant y$, so, $u[\bar{w}|_i(y-x)]\leqslant
u[\bar{w}|_i\,y]$ and
\[
 u[\bar{w}|_i(y-x)]\leqslant u[\bar{w}|_i\,x]\curlyvee u[\bar{w}|_i\,y].
\]
This and \eqref{e-58} guarantee the existence of an element
$$
(u[\bar{w}|_i(x\curlyvee y)]-u[\bar{w}|_i(y-x)])\curlyvee
u[\bar{w}|_i(y-x)]
$$
such that
\begin{equation}\label{e-59}
 (u[\bar{w}|_i(x\curlyvee y)]-u[\bar{w}|_i(y-x)])\curlyvee
u[\bar{w}|_i(y-x)]\leqslant u[\bar{w}|_i\,x]\curlyvee
u[\bar{w}|_i\,y].
\end{equation}
Since $u[\bar{w}|_i(y-x)]\leqslant u[\bar{w}|_i\,y]\leqslant
u[\bar{w}|_i(x\curlyvee y)]$, the last inequality and \eqref{e-51}
imply
$$
(u[\bar{w}|_i(x\curlyvee y)]-u[\bar{w}|_i(y-x)])\curlyvee
u[\bar{w}|_i(y-x)]=u[\bar{w}|_i(x\curlyvee y)],
$$
which together with \eqref{e-59} gives
\[
u[\bar{w}|_i(x\curlyvee y)]\leqslant u[\bar{w}|_i\,x]\curlyvee
u[\bar{w}|_i\,y].
\]
Comparing this inequality with \eqref{e-57} we obtain
\eqref{e-56}.
\end{proof}
\begin{collolary}\label{C-2}
If for some $x,y\in G$ an element $x\curlyvee y$ exists, then for
any polynomial $t\in T_n(G)$ an element $t(x)\curlyvee t(y)$ also
exists and $t(x\curlyvee y)=t(x)\curlyvee t(y)$.
\end{collolary}
\begin{proposition}\label{P-8}
For all $x,y\in G$ and all polynomials $t_1,t_2\in T_n(G)$ we
have:
$$
t_1(x\curlywedge y)\curlywedge t_2(x-y)=0.
$$
\end{proposition}
\begin{proof}
Let $t_1(x\curlywedge y)\curlywedge t_2(x-y)=h$. Obviously
$h\leqslant t_1(x\curlywedge y)$ and $h\leqslant t_2(x-y)$. Since
$t_2(x-y)\leqslant t_2(x)$, we have $h\leqslant t_2(x)$. Thus,
$x\curlywedge y\leqslant x$, $h\leqslant t_1(x\curlywedge y)$ and
$h\leqslant t_2(x)$. This, in view of Proposition \ref{P-3} and
\eqref{e-13}, gives $h\leqslant t_2(x\curlywedge y)$.
Consequently,
\begin{equation}\label{e-60}
  h\leqslant t_2(x-y)\curlywedge t_2(x\curlywedge y).
\end{equation}
Further,
\begin{eqnarray}
 t_2(x-y)-t_2(x\curlywedge
y)&\stackrel{\eqref{e-39}}{=}&(t_2(x)-t_2(x\curlywedge
y))-t_2(x\curlywedge y)\nonumber\\
&\stackrel{\eqref{e-20}}{=}&t_2(x)-t_2(x\curlywedge
y)\stackrel{\eqref{e-39}}{=}t_2(x-y).\nonumber
\end{eqnarray}
Therefore,
$$t_2(x-y)\curlywedge t_2(x\curlywedge
y)\stackrel{\eqref{e-25}}{=}t_2(x-y)-(t_2(x-y)-t_2(x\curlywedge
y))=t_2(x-y)-t_2(x-y)=0,
$$
which together with \eqref{e-60} implies $h\leqslant 0$. Hence
$h=0$. This completes the proof.
\end{proof}
\begin{proposition}\label{P-9}
For all $x,y,z,g\in G$ and all polynomials $t_1,t_2\in T_n(G)$ the
following conditions are valid:
\begin{eqnarray}
&&\label{e-61} t_1(x\curlywedge y)\curlywedge
t_2(y)=t_1(x\curlywedge y)\curlywedge t_2(x\curlywedge y),\\
&&\label{e-62} t_1(x\curlywedge y\curlywedge z)\curlywedge
t_2(y)\leqslant t_1(x\curlywedge y)\curlywedge t_2(y\curlywedge z),\\
&&\label{e-63} g\leqslant t_1(x\curlywedge y)\wedge g\leqslant
t_2(y\curlywedge z)\longrightarrow g\leqslant t_2(x\curlywedge
y\curlywedge z).
\end{eqnarray}
\end{proposition}
\begin{proof}
To prove \eqref{e-61} observe first that for $z=t_1(x\curlywedge
y)\curlywedge t_2(y)$ we have $z\leqslant t_1(x\curlywedge y)$ and
$z\leqslant t_2(y)$. Since the relation $\omega$ is weakly steady
and $x\curlywedge y\leqslant y$, from the above we conclude
$z\leqslant t_2(x\curlywedge y)$, i.e., $t_1(x\curlywedge
y)\curlywedge t_2(y)\leqslant t_2(x\curlywedge y)$. This, by
\eqref{e-31}, implies $t_1(x\curlywedge y)\curlywedge
t_2(y)\leqslant t_1(x\curlywedge y)\curlywedge t_2(x\curlywedge
y)$.

On the other side, the stability of $\omega$ and $x\curlywedge
y\leqslant y$ imply $t_2(x\curlywedge y)\leqslant t_2(y)$ for
every $t_2\in T_n(G)$. Hence, $t_1(x\curlywedge y)\curlywedge
t_2(x\curlywedge y)\leqslant t_1(x\curlywedge y)\curlywedge
t_2(y)$ by \eqref{e-31}. This completes the proof of \eqref{e-61}.

Further: $t_1(x\curlywedge y\curlywedge z)\curlywedge
t_2(y)=t_1((x\curlywedge z)\curlywedge y)\curlywedge
t_2(y)\stackrel{\eqref{e-61}}{=}t_1((x\curlywedge z)\curlywedge
y)\curlywedge t_2((x\curlywedge z)\curlywedge y)\leqslant
t_1(x\curlywedge y)\curlywedge t_2(y\curlywedge z)$ proves
\eqref{e-62}.

Finally, let $g\leqslant t_1(x\curlywedge y)$ and $g\leqslant
t_2(y\curlywedge z)$. Then
\begin{eqnarray}
&& g\leqslant t_1(x\curlywedge y)\curlywedge t_2(y\curlywedge
z)\stackrel{\eqref{e-28}}{=}t_1(x\curlywedge y)\curlywedge
(t_2(y)-t_2(y-z))\nonumber\\&&
\stackrel{\eqref{e-34}}{=}\Big(t_1(x\curlywedge y)\curlywedge
t_2(y)\Big)-\Big(t_1(x\curlywedge y)\curlywedge
t_2(y-z)\Big)\nonumber\\&&
\stackrel{\eqref{e-61}}{=}\Big(t_1(x\curlywedge y)\curlywedge
t_2(x\curlywedge y)\Big)-\Big(t_1(x\curlywedge y)\curlywedge
t_2(y-z)\Big)\nonumber\\&&
\stackrel{\eqref{e-34}}{=}t_1(x\curlywedge
y)\curlywedge\Big(t_2(x\curlywedge y)- t_2(y-z)\Big)\leqslant
t_2(x\curlywedge y)- t_2(y-z)\nonumber\\&&
\stackrel{\eqref{e-40}}{\leqslant}t_2\Big((x\curlywedge
y)-(y-z)\Big)\stackrel{\eqref{e-36}}{=}t_2(x\curlywedge
y\curlywedge z ).\nonumber
\end{eqnarray}
This proves \eqref{e-63} and completes the proof of our
proposition.
\end{proof}

\begin{collolary}\label{C-3}
For all $x,y,z\in G$ and all polynomials $t_1,t_2\in T_n(G)$ we
have:
\begin{equation}\label{e-64}
t_1(x\curlywedge y\curlywedge
z)\curlywedge t_2(y)= t_1(x\curlywedge y)\curlywedge
t_2(y\curlywedge z).
\end{equation}
\end{collolary}
\begin{proof}
We have $t_1(x\curlywedge y)\curlywedge t_2(y\curlywedge
z)\leqslant t_1(x\curlywedge y)$ and $t_1(x\curlywedge
y)\curlywedge t_2(y\curlywedge z)\leqslant t_2(y\curlywedge z)$,
so by \eqref{e-63} we obtain $t_1(x\curlywedge y)\curlywedge
t_2(y\curlywedge z)\leqslant t_1(x\curlywedge y\curlywedge z)$.
Considering now that $t_1(x\curlywedge y)\curlywedge
t_2(y\curlywedge z)\leqslant t_2(y\curlywedge z)\leqslant t_2(y)$,
by \eqref{e-30}, we get $t_1(x\curlywedge y)\curlywedge
t_2(y\curlywedge z)\leqslant t_1(x\curlywedge y\curlywedge
z)\curlywedge t_2(y)$. Taking now into account the condition
\eqref{e-62} we obtain \eqref{e-64}.
\end{proof}

\medskip {\bf 5.} Let $(G,o,-,0)$ be a subtraction Menger algebra of rank $n$.

\begin{definition}\rm
By a {\it determining pair} of a subtraction Menger algebra
$(G,o,-,0)$ of rank $n$ we mean an ordered pair
$(\varepsilon^*,W)$, where $\varepsilon$ is a $v$-regular
equivalence relation defined on $(G,o)$,
$\varepsilon^*=\varepsilon\cup\{(e_1,e_1),\ldots,(e_n,e_n)\}$,
$e_1,\ldots,e_n$ are selectors of a unitary extension $(G^*,o^*)$
of $(G,o)$ and $W$ is the empty set or an $l$-ideal of $(G,o)$
which is an $\varepsilon$-class.
\end{definition}
\begin{definition}\rm
A non-empty subset $F$ of a subtraction Menger algebra $(G,o,-,0)$
of rank $n$ is called a {\it filter } if:

\smallskip $
\begin{array}{rl}
1)& 0\not\in F;\\
2)& x\in F\wedge x\leqslant y\longrightarrow y\in F;\\
3)& x\in F\wedge y\in F\longrightarrow x\curlywedge y\in F
\end{array}$

\noindent for all $x,y\in G$.
\end{definition}

If $a,b\in G$ and $a\nleqslant b$, then $[\,a)=\{x\in
G\,|\,a\leqslant x \}$ is a filter with $a\in[\,a)$ and
$b\not\in[\,a)$. By Zorn's Lemma the collection of filters which
contain an element $a$, but do not contain an element $b$, has a
maximal element which is denoted by $F_{a,b}$. Using this filter
we define the following three sets:
\begin{eqnarray}
&& W_{a,b}=\{x\in G\,|\,(\forall t\in T_n(G))\,t(x)\not\in
F_{a,b}\},\nonumber\\
&& \varepsilon_{a,b}=\{(x,y)\in G\times
G\,|\,x\curlywedge y\not\in W_{a,b}\vee x,y\in W_{a,b}\},\nonumber\\
&&\varepsilon^*_{a,b}=\varepsilon_{a,b}\cup\{(e_1,e_1),\ldots,(e_n,e_n)\}.\nonumber
\end{eqnarray}
\begin{proposition}\label{P-10}
For any $a,b\in G$, the pair $(\varepsilon^*_{a,b},W_{a,b})$ is
the determining pair of the algebra $(G,o,-,0)$.
\end{proposition}
\begin{proof}
First we show that $\varepsilon_{a,b}$ is an equivalence relation
on $G$. It is clear that this relation is reflexive and symmetric.
To prove its transitivity let $(x,y),(y,z)\in\varepsilon_{a,b}$.
We have four possibilities:

$(a)$ \ $x\curlywedge y\not\in W_{a,b}\wedge\,y\curlywedge
z\not\in W_{a,b}$,

$(b)$ \ $x\curlywedge y\not\in W_{a,b}\wedge\,y,z\in W_{a,b}$,

$(c)$ \ $x,y\in W_{a,b}\wedge\,y\curlywedge z\not\in W_{a,b}$,

$(d)$ \ $x,y\in W_{a,b}\wedge\,y,z\in W_{a,b}$.

In the case $(a)$ we have $t_1(x\curlywedge y),t_2(y\curlywedge
z)\in F_{a,b}$ for some $t_1,t_2\in T_n(G)$. Since $F_{a,b}$ is a
filter, then, obviously, $t_1(x\curlywedge y)\curlywedge
t_2(y\curlywedge z)\in F_{a,b}$. This, according to \eqref{e-64},
implies $t_1(x\curlywedge y\curlywedge z)\curlywedge t_2(y)\in
F_{a,b}$. But $t_1(x\curlywedge y\curlywedge z)\curlywedge
t_2(y)\leqslant t_1(x\curlywedge z)$, hence also $t_1(x\curlywedge
z)\in F_{a,b}$, i.e., $x\curlywedge z\not\in W_{a,b}$. Thus,
$(x,z)\in\varepsilon_{a,b}$.

In the case $(b)$ from $x\curlywedge y\not\in W_{a,b}$ it follows
$t(x\curlywedge y)\in F_{a,b}$ for some polynomial $t\in T_n(G)$.
But $x\curlywedge y\leqslant y$, and consequently $t(x\curlywedge
y)\leqslant t(y)$. Thus $t(y)\in F_{a,b}$, i.e., $y\not\in
W_{a,b}$, which is a contradiction. Hence the case $(b)$ is
impossible. Analogously we can show that also the case $(c)$ is
impossible. The case $(d)$ is obvious, because in this case
$x,z\in W_{a,b}$ which means that $(x,z)\in\varepsilon_{a,b}$.
This completes the proof that $\varepsilon_{a,b}$ is transitive.

Moreover, if $x\in W_{a,b}$, then $t(x)\not\in F_{a,b}$ for every
$t\in T_n(G)$. In particular, for all
$t(x)=t'(u[\bar{w}|_i\,x])\in T_n(G)$ we have
$t'(u[\bar{w}|_i\,x])\not\in F_{a,b}$. Thus, $u[\bar{w}|_i\,x]\in
W_{a,b}$ for every $i=1,\ldots,n$. Hence, $W_{a,b}$ is an
$i$-ideal of $(G,o)$, and consequently, an $l$-ideal. It is clear
that $W_{a,b}$ is an $\varepsilon_{a,b}$-class.

Next, we prove that the relation $\varepsilon_{a,b}$ is
$v$-regular. Let $x\equiv y(\varepsilon_{a,b})$. Then
$x\curlywedge y\not\in W_{a,b}$ or $x,y\in W_{a,b}$. In the case
$x,y\in W_{a,b}$ we obtain $u[\bar{w}|_i\,x],u[\bar{w}|_i\,y]\in
W_{a,b}$ because $W_{a,b}$ is an $l$-ideal of $(G,o)$. Thus,
$u[\bar{w}|_i\,x]\equiv u[\bar{w}|_i\,y](\varepsilon_{a,b})$. In
the case $x\curlywedge y\not\in W_{a,b}$ elements
$u[\bar{w}|_i\,x]$, $u[\bar{w}|_i\,y]$ belong or not belong to
$W_{a,b}$ simultaneously. Indeed, if $u[\bar{w}|_i\,x]$,
$u[\bar{w}|_i\,y]\in W_{a,b}$, then obviously
$u[\bar{w}|_i\,x]\equiv u[\bar{w}|_i\,y](\varepsilon_{a,b})$. Now,
if $u[\bar{w}|_i\,x]\not\in W_{a,b}$, then $t(u[\bar{w}|_i\,x])\in
F_{a,b}$ for some $t\in T_n (G)$. Since $x\curlywedge y\not\in
W_{a,b}$, then also $t_1(x\curlywedge y)\in F_{a,b}$ for some
$t_1\in T_n(G)$. Thus $t_1(x\curlywedge y) \curlywedge
t(u[\bar{w}|_i\,x])\in F_{a,b}$, which, by \eqref{e-61}, implies
$t_1(x\curlywedge y)\curlywedge t(u[\bar{w}|_i(x\curlywedge
y)])\in F_{a,b}$. But $t_1(x\curlywedge y)\curlywedge
t(u[\bar{w}|_i(x\curlywedge y)]) \leqslant t(u[\bar{w}|_i\,y])$,
hence $t(u[\bar{w}|_i\,y])\in F_{a,b}$, i.e.,
$u[\bar{w}|_i\,y]\not\in W_{a,b}$. So, we have shown that
$x\curlywedge y\not\in W_{a,b}$ and $u[\bar{w}|_i\,x]\not\in
W_{a,b}$ imply $u[\bar{w}|_i\,y]\not\in W_{a,b}$. Similarly we can
show that $x\curlywedge y\not\in W_{a,b}$ and
$u[\bar{w}|_i\,y]\not\in W_{a,b}$ imply $u[\bar{w}|_i\,x]\not\in
W_{a,b}$. Therefore, we have proved that in the case $x\curlywedge
y\not\in W_{a,b}$ elements $u[\bar{w}|_i\,x]$, $u[\bar{w}|_i\,y]$
belong or not belong to $W_{a,b}$ simultaneously.

So, if for $x\curlywedge y\not\in W_{a,b}$ we have
$u[\bar{w}|_i\,x], u[\bar{w}|_i\,y]\in W_{a,b}$, then clearly
$u[\bar{w}|_i\,x]\equiv u[\bar{w}|_i\,y](\varepsilon_{a,b})$.
Therefore assume that $u[\bar{w}|_i\,x]\not\in W_{a,b}$ (hence
$u[\bar{w}|_i\,y]\not\in W_{a,b}$). Thus, $x\curlywedge y\not\in
W_{a,b}$, $u[\bar{w}|_i\,x]\not\in W_{a,b}$, i.e., $t(x\curlywedge
y)\in F_{a,b}$, $t_1(u[\bar{w}|_i\,x])\in F_{a,b}$ for some
$t,t_1\in T_n(G)$. Hence, $t(y\curlywedge x\curlywedge
y)\curlywedge t_1(u[\bar{w}|_i\,x])\in F_{a,b}$. From this,
according to \eqref{e-64}, we obtain $t(y\curlywedge x)\curlywedge
t_1(u[\bar{w}|_i(x\curlywedge y)])\in F_{a,b}$. This implies
$t_1(u[\bar{w}|_i(x\curlywedge y)])\in F_{a,b}$. Since
$u[\bar{w}|_i(x\curlywedge y)]\leqslant u[\bar{w}|_i\,x]$ and
$u[\bar{w}|_i(x\curlywedge y)]\leqslant u[\bar{w}|_i\,y]$, we have
$u[\bar{w}|_i(x\curlywedge y)]\leqslant
u[\bar{w}|_i\,x]\curlywedge u[\bar{w}|_i\,y]$, which, by the
stability of $\omega$ gives $t_1(u[\bar{w}|_i(x\curlywedge y)])
\leqslant t_1(u[\bar{w}|_i\,x]\curlywedge u[\bar{w}|_i\,y])$.
Consequently, $t_1(u[\bar{w}|_i\,x]\curlywedge
u[\bar{w}|_i\,y])\in F_{a,b}$, so $u[\bar{w}|_i\,x]\curlywedge
u[\bar{w}|_i\,y]\not\in W_{a,b}$, i.e., $u[\bar{w}|_i\,x]\equiv
u[\bar{w}|_i\,y](\varepsilon_{a,b})$. In this way we have proved
that the relation $\varepsilon_{a,b}$ is $i$-regular for every
$i=1,\ldots,n$. Thus it is $v$-regular.
\end{proof}

\begin{proposition}\label{P-11} All equivalence classes of $\varepsilon_{a,b}$,
except of $W_{a,b}$, are filters.
\end{proposition}
\begin{proof}
Indeed, let $H\ne W_{a,b}$ be an arbitrary class of
$\varepsilon_{a,b}$. If $x\in H$ and $x\leqslant y$, then
$x\curlywedge y=x\not\in W_{a,b}$, consequently,
$(x,y)\in\varepsilon_{a,b}$. Hence, $y\in H$. Further, let $x,y\in
H$, then $(x,y)\in\varepsilon_{a,b}$. Thus $x\curlywedge y\not\in
W_{a,b}$, i.e., $t(x\curlywedge y)\in F_{a,b}$ for some $t\in
T_n(G)$. But $x\curlywedge y=x\curlywedge(x\curlywedge y)$, hence,
$t(x\curlywedge(x\curlywedge y))\in F_{a,b}$ and $x\curlywedge
(x\curlywedge y)\not\in W_{a,b}$. So $x\equiv x\curlywedge
y(\varepsilon_{a,b})$. This implies $x\curlywedge y\in H$. Thus,
we have shown that $H$ is a filter.
\end{proof}

\begin{proposition}\label{P-12}
If $x\curlyvee y$ exists for some $x,y\in W_{a,b}$, then
$x\curlyvee y\in W_{a,b}$.
\end{proposition}
\begin{proof}
Let $x\curlyvee y$ exists for some $x,y\in W_{a,b}$. If
$x\curlyvee y\not\in W_{a,b}$, then $t(x\curlyvee y)\in F_{a,b}$
for some $t\in T_n(G)$, and, according to Corollary \ref{C-2},
$t(x\curlyvee y)=t(x)\curlyvee t(y)$. If $t(x)\not\in F_{a,b}$,
then $F_{a,b}$ is a proper subset of the set
$$
U =\{u\in G\,|\,(\exists z\in F_{a,b})\,z\curlywedge t(x)\leqslant
u\}
$$
because $t(x)\in U$.

We show that $U$ is a filter. $0\not\in U$ because, by
\eqref{e-15}, we have $0\leqslant z\curlywedge t(x)$ for any $z\in
F_{a,b}$. Let $s\in U$ and $s\leqslant r$. Then $z\curlywedge
t(x)\leqslant s$ for some $z\in F_{a,b}$. Consequently,
$z\curlywedge t(x)\leqslant r$, so $r\in U$. Now let $s\in U$ and
$r\in U$, i.e., $z_1\curlywedge t(x)\leqslant s$ and
$z_2\curlywedge t(x)\leqslant r$ for some $z_1,z_2\in F_{a,b}$.
Since $F_{a,b}$ is a filter, we have $z_1\curlywedge z_2\in
F_{a,b}$. Hence, $(z_1\curlywedge z_2)\curlywedge t(x)\leqslant
s\curlywedge r$, which implies $s\curlywedge r\in U$. Thus $U$ is
a filter. But by assumption $F_{a,b}\subset U$ is a maximal
filter, which does not contain $b$, so $b\in U$. Consequently,
$z_1\curlywedge t(x)\leqslant b$ for some $z_1\in F_{a,b}$.
Similarly, if $t(y)\not\in F_{a,b}$, then $z_2\curlywedge
t(y)\leqslant b$ for some $z_2\in F_{a,b}$. This implies
$z\curlywedge t(x)\leqslant b$ and $z\curlywedge t(y)\leqslant b$
for $z=z_1\curlywedge z_2$. Hence $(z\curlywedge t(x))\curlyvee
(z\curlywedge t(y))$ exists and
\[
(z\curlywedge t(x))\curlyvee(z\curlywedge
t(y))=z\curlywedge(t(x)\curlyvee t(y)) =z\curlywedge t(x\curlyvee
y)\in F_{a,b}
\]
by \eqref{e-47}. But by \eqref{e-50} we have $(z\curlywedge
t(x))\curlyvee (z\curlywedge t(y))\leqslant b$, so $z\curlywedge
t(x\curlyvee y) \leqslant b$. Since $z\curlywedge t(x\curlyvee
y)\in F_{a,b}$, then, obviously, $b\in F_{a,b}$, which is
impossible. So, $t(x)\in F_{a,b}$ or $t(y)\in F_{a,b}$, hence
$x\not\in W_{a,b}$ or $y\not\in W_{a,b}$, which is contrary to the
assumption that $x,y\in W_{a,b}$. Thus, the assumption that
$x\curlyvee y\not\in W_{a,b}$ is incorrect. Therefore $x\curlyvee
y\in W_{a,b}$.
\end{proof}

\medskip {\bf 6.}
Each homomorphism of a Menger algebra $(G,o)$ of rank $n$ into a
Menger algebra $(\mathcal{F}(A^n,A),\mathrm{O})$ is called a {\it
representation  by $n$-place functions}. Thus, $P:G\to
\mathcal{F}(A^n,A)$ is a representation, if
$$
P(x[y_1\ldots y_n])=P(x)[P(y_1)\ldots P(y_n)]
$$
for all $x,y_1,\ldots,y_n\in G$. A representation which is an
isomorphism is called {\it faithful} (cf. \cite{Dudtro2,Dudtro3,
SchTro}). A representation $P$ of $(G,o)$ is a representation of
$(G,o,-,0)$ if
$$
P(x-y)=P(x)\setminus P(y) \ \ \ {\rm and } \ \ \ P(0)=\varnothing
$$
for all $x,y\in G$.

Let $(P_i)_{i\in I}$ be the family of representations of a
subtraction Menger algebra $(G,o,-,0)$ of rank $n$ by $n$-place
functions defined on pairwise disjoint sets $(A_i)_{i\in I}$. By
the {\it sum} of the family $(P_i)_{i\in I}$ we mean the map
$P\colon g\mapsto P(g)$, denoted by $\sum_{i\in I}P_i$, where
$P(g)$ is an $n$-place function on $A=\bigcup_{i\in I}A_i$ defined
by $P(g)= \bigcup_{i\in I}P_i(g)$. It is clear (cf.
\cite{Dudtro2,Dudtro3}) that $P$ is a representation of
$(G,o,-,0)$.

Similarly as in \cite{Dudtro2,Dudtro3} with each determining pair
$(\varepsilon^*,W)$ we can associate the so-called {\it simplest
representation} $P_{(\varepsilon^*,W)}$ of $(G,o)$ which assigns
to each element $g\in G$ an $n$-place function
$P_{(\varepsilon^*,W)}(g)$ defined on
$\mathcal{H}=\mathcal{H}_0\cup \{\{e_1\},\ldots,\{e_n\}\}$, where
$\mathcal{H}_0$ is the set of all $\varepsilon$-classes of $G$
different from $W$ such that
\[
(H_1,\ldots,H_n,H)\in P_{(\varepsilon,W)}(g)\longleftrightarrow
g[H_1\ldots H_n]\subset H,
\]
for $(H_1,\ldots,H_n)\in\mathcal{H}^n_0\cup
\{(\{e_1\},\ldots,\{e_n\})\}$ and $H\in\mathcal{H}$.

\begin{theorem}\label{T-2}
Each subtraction Menger algebra of rank $n$ is isomorphic to some
difference Menger algebra of $n$-place functions.
\end{theorem}
\begin{proof}
Let $(G,o,-,0)$ be a subtraction Menger algebra of rank $n$. Then
the sum
\[
P=\sum_{a,b\in G,\,a\nleqslant b}P_{(\varepsilon^*_{a,b},W_{a,b})}
\]
of the family $\Big(P_{(\varepsilon^*_{a,b},W_{a,b})}\Big)_{a,b\in
G,\,a\nleqslant b}$ of simplest representations of $(G,o)$ is a
representation of $(G,o)$.

Now we show that $P$ is a representation of $(G,o,-,0)$. Let
$\mathcal{H}_0$ be the set of all $\varepsilon_{a,b}$-classes of
$G$ different from $W_{a,b}$. Consider $H_1,\ldots,H_n,H\in
\mathcal{H}$, where $\mathcal{H}=\mathcal{H}_0\cup
\{\{e_1\},\ldots,\{e_n\}\}$, such that $(H_1,\ldots,H_n,H)\in
P_{(\varepsilon^*_{a,b},W_{a,b})}(g_1-g_2)$ for some $g_1,g_2\in
G$. Then, obviously, $(g_1-g_2)[H_1\ldots H_n]\subset H\neq
W_{a,b}$. Thus $(g_1-g_2)[\bar{x}]\in H$ for each $\bar{x}\in
H_1\times\cdots\times H_n$, which, by \eqref{e-11}, gives
$g_1[\bar{x}]-g_2[\bar{x}]\in H$. But
$g_1[\bar{x}]-g_2[\bar{x}]\leqslant g_1[\bar{x}]$ and $H$ is a
filter (Proposition \ref{P-11}), hence $g_1[\bar{x}]\in H$. Thus
$(g_1[\bar{x}]-g_2[\bar{x}])\curlywedge g_2[\bar{x}]=0$, by
\eqref{e-33}. Consequently,
$(g_1[\bar{x}]-g_2[\bar{x}])\curlywedge g_2[\bar{x}]\in W_{a,b}$,
because the other $\varepsilon_{a,b}$-classes as filters do not
contain $0$. This means that $g_1[\bar{x}]-g_2[\bar{x}]\not\equiv
g_2[\bar{x}](\varepsilon_{a,b})$. Hence, $g_2[\bar{x}]\not\in H$.
Therefore $g_1[H_1\ldots H_n]\subset H$ and $g_2[h_1\ldots
H_n]\cap H=\varnothing$, which implies
$$
(H_1,\ldots,H_n,H)\in
P_{(\varepsilon^*_{a,b},W_{a,b})}(g_1)\setminus
P_{(\varepsilon^*_{a,b},W_{a,b})}(g_2).
$$
In this way, we have proved the inclusion
\begin{equation}\label{e-65}
P_{(\varepsilon^*_{a,b},W_{a,b})}(g_1-g_2)\subset
P_{(\varepsilon^*_{a,b},W_{a,b})}(g_1)\setminus
P_{(\varepsilon^*_{a,b},W_{a,b})}(g_2).
\end{equation}

To show the reverse inclusion let
$$
(H_1,\ldots,H_n,H)\in
P_{(\varepsilon^*_{a,b},W_{a,b})}(g_1)\setminus
P_{(\varepsilon^*_{a,b},W_{a,b})}(g_2).
$$
Then $(H_1,\ldots,H_n,H)\in
P_{(\varepsilon^*_{a,b},W_{a,b})}(g_1)$ and
$(H_1,\ldots,H_n,H)\not\in
P_{(\varepsilon^*_{a,b},W_{a,b})}(g_2)$, i.e., $g_1[H_1\ldots
H_n]\subset H$ and $g_2[H_1\ldots H_n]\cap H=\varnothing$. Thus
$g_1[\bar{x}]\in H$ and $g_2[\bar{x}]\not\in H$ for all
$\bar{x}\in H_1\times\cdots\times H_n$. Since from
$g_1[\bar{x}]\curlywedge g_2[\bar{x}]\not\in W_{a,b}$, it follows
$g_1[\bar{x}]\equiv g_2[\bar{x}](\varepsilon_{a,b})$ and
$g_2[\bar{x}]\in H$, which is a contradiction, we conclude that
$g_1[\bar{x}]\curlywedge g_2[\bar{x}]\in W_{a,b}$.

If $g_1[\bar{x}]-g_2[\bar{x}]\in W_{a,b}$, then, by \eqref{e-53}
and Proposition \ref{P-12}, we obtain
$g_1[\bar{x}]=(g_1[\bar{x}]\curlywedge
g_2[\bar{x}])\curlyvee(g_1[\bar{x}]-g_2[\bar{x}])\in W_{a,b}$.
Consequently, $g_1[\bar{x}]\in W_{a,b}$, which is impossible
because $g_1[\bar{x}]\in H$. Thus,
$(g_1[\bar{x}]-g_2[\bar{x}])\curlywedge
g_1[\bar{x}]=g_1[\bar{x}]-g_2[\bar{x}]\not\in W_{a,b}$. Hence,
$g_1[\bar{x}]-g_2[\bar{x}]\equiv g_1[\bar{x}](\varepsilon_{a,b})$.
This implies $(g_1-g_2)[\bar{x}]=g_1[\bar{x}]-g_2[\bar{x}]\in H$.
Therefore, $(g_1-g_2)[H_1\ldots H_n]\subset H$, i.e.,
$(H_1,\ldots,H_n,H)\in
P_{(\varepsilon^*_{a,b},W_{a,b})}(g_1-g_2)$. So, we have proved
\[
P_{(\varepsilon^*_{a,b},W_{a,b})}(g_1)\setminus
P_{(\varepsilon^*_{a,b},W_{a,b})}(g_2)\subset
P_{(\varepsilon^*_{a,b},W_{a,b})}(g_1-g_2).
\]
This together with \eqref{e-65} proves
\[
P_{(\varepsilon^*_{a,b},W_{a,b})}(g_1-g_2)=
P_{(\varepsilon^*_{a,b},W_{a,b})}(g_1)\setminus
P_{(\varepsilon^*_{a,b},W_{a,b})}(g_2),
\]
which means that $P(g_1-g_2)=P(g_1)\setminus P(g_2)$ for
$g_1,g_2\in G$. Further, $P(0)=P(0-0)=P(0)\setminus
P(0)=\varnothing$. So, $P$ is a representation of $(G,o,-,0)$ by
$n$-place functions.

We show that this representation is faithful. Let $P(g_1)=P(g_2)$
for some $g_1,g_2\in G$. If $g_1\neq g_2$, then both inequalities
$g_1\leqslant g_2$ and $g_2\leqslant g_1$ at the same time are
impossible. Suppose that $g_1\nleqslant g_2$. Then $g_1\in
F_{g_1,g_2}$ and, consequently,
$$
(\{e_1\},\ldots,\{e_n\},F_{g_1,g_2})\in
P_{(\varepsilon^*_{g_1,g_2},W_{g_1,g_2}})(g_2).
$$
Since
$P_{(\varepsilon^*_{g_1,g_2},W_{g_1,g_2}})(g_1)=P_{(\varepsilon^*_{g_1,g_2},W_{g_1,g_2}})(g_2)$,
then, obviously,
$$
(\{e_1\},\ldots,\{e_n\},F_{g_1,g_2})\in
P_{(\varepsilon^*_{g_1,g_2},W_{g_1,g_2}})(g_2).
$$
Thus $\{g_2\}=g_2[\{e_1\}\ldots\{e_n\}]\subset F_{g_1,g_2}$, hence
$g_2\in F_{g_1,g_2}$. This is a contradiction because
$F_{g_1,g_2}$ is a filter containing $g_1$ but not containing
$g_2$. The case $g_2\nleqslant g_1$ is analogous. So, the
supposition $g_1\ne g_2$ is not true. Hence $g_1=g_2$ and $P$ is a
faithful representation. The theorem is proved.
\end{proof}

\begin{minipage} {60mm}
\begin{flushleft}
Dudek~W. A.\\ Institute of Mathematics and Computer Science\\
Wroclaw University of Technology\\ 50-370 Wroclaw\\ Poland\\
Email: dudek@im.pwr.wroc.pl
\end{flushleft}
\end{minipage}
\hfill
\begin{minipage} {60mm}
\begin{flushleft}
 Trokhimenko~V. S.\\ Department of Mathematics\\
 Pedagogical University\\
 21100 Vinnitsa\\
 Ukraine\\
 Email: vtrokhim@sovamua.com
\end{flushleft}
\end{minipage}
\end{document}